\documentclass[12pt,oneside]{amsart}
 \usepackage{mathcomp,longtable, hhline}
\usepackage{amsfonts,amssymb,mathrsfs}
\usepackage{url}
 \usepackage[matrix,arrow,curve]{xy}

\makeatletter \@addtoreset{equation}{section} \makeatother

\newcommand{\Pic}{\operatorname{Pic}}
\newcommand{\Spec}{\operatorname{Spec}}
\newcommand{\Diff}{\operatorname{Diff}}
\newcommand{\Sing}{\operatorname{Sing}}
\newcommand{\Supp}{\operatorname{Supp}}

\newcommand{\codim}{\operatorname{codim}}
\newcommand{\Bs}{\operatorname{Bs}}
\newcommand{\Fix}{\operatorname{Fix}}
\newcommand{\F}{\operatorname{F}}
\newcommand{\Cr}{\operatorname{Cr}}
\newcommand{\Cl}{\operatorname{Cl}}
\newcommand{\Aut}{\operatorname{Aut}}
\newcommand{\Bir}{\operatorname{Bir}}

\newcommand{\dd}{\operatorname{d}}
\newcommand{\g}{\operatorname{g}}

\newcommand{\Proj}{\operatorname{Proj}}

\newcommand{\Gal}{\operatorname{Gal}}

\newcommand{\A}{{\mathbb A}}

\newcommand{\CC}{\mathbb{C}}
\newcommand{\LL}{\mathbb{L}}
\newcommand{\FF}{\mathbb{F}}
\newcommand{\QQ}{\mathbb{Q}}

\newcommand{\ZZ}{\mathbb{Z}}

\newcommand{\PP}{\mathbb{P}}
\newcommand{\bA}{\mathbb{A}}
\newcommand{\OOO}{{\mathscr{O}}} 
 
\newcommand{\MMM}{{\mathscr{M}}} 
\newcommand{\LLL}{{\mathscr{L}}} 

\newcommand{\EEE}{{\mathscr{E}}}
\renewcommand{\emptyset}{\varnothing}

\newcommand{\type}[1]{$\mathrm{(#1)}$}

\newcommand{\qq}{\mathbin{\sim_{\scriptscriptstyle{\QQ}}}}

\newcounter{NO}

\renewcommand\labelenumi{(\roman{enumi})}
\renewcommand\theenumi{(\roman{enumi})}

\newcommand{\comment}[1]{}

\newcommand{\xref}[1]{{\rm \ref{#1}}}

\newtheorem{theorem}{Theorem}

\numberwithin{theorem}{section}
\numberwithin{equation}{theorem}

\newtheorem{mtheorem}[theorem]{}
\newtheorem{stheorem}[equation]{}

\theoremstyle{definition}
\newtheorem{case}[theorem]{}
\newtheorem{scase}[equation]{}

\newcounter{NN}\numberwithin{NN}{section}

\begin{document}

\title{On birational involutions of $\PP^3$}

\author{Yuri Prokhorov}
\thanks{
I acknowledge partial supports by
RFBR grants No. 11-01-00336-a, 11-01-92613-KO\_a, the grant of
Leading Scientific Schools No. 4713.2010.1, Simons-IUM fellowship,
 and 
AG Laboratory SU-HSE, RF government 
grant ag. 11.G34.31.0023.
}

\address{
Department 
of Algebra, Faculty of Mathematics, Moscow State
University, Moscow 117234, Russia
\newline\indent
Laboratory of Algebraic Geometry, SU-HSE, 
7 Vavilova Str., Moscow 117312, Russia
}
\email{prokhoro@gmail.com}

\begin{abstract}
Let $X$ be a rationally connected three-dimensional 
algebraic variety and let $\tau$
be an element of order two in the group of its birational selfmaps. 
Suppose that there exists a non-uniruled divisorial component of the $\tau$-fixed point locus. 
Using the equivariant minimal model program we give a 
rough classification of such elements. 
\end{abstract}

\maketitle

\section{Introduction}\label{sect-intr}
Let $\Bbbk$ be an algebraically closed field of characteristic $0$.
The \textit{Cremona group} $\Cr_n(\Bbbk)$ is the group of birational transformations of 
$\PP^n_{\Bbbk}$, or equivalently the group of $\Bbbk$-automorphisms of the field $\Bbbk(x_1,\dots,x_n)$. 
For $n=2$ the Cremona group and its subgroups were studied extensively 
in the classical literature as well as in many recent works.
In particular, the classification of finite subgroups in $\Cr_2(\Bbbk)$ was started in the works of 
E. Bertini and was almost completed recently by I. Dolgachev and V. Iskovskikh \cite{Dolgachev-Iskovskikh}.
For the Cremona group in three variables the situation is much more complicated.
There are only a few partial results in this direction, see e.g. \cite{Prokhorov2011a}, \cite{Prokhorov2009e}.

The modern approach to the classification of finite subgroups in Cremona groups
is based on the simple observation that any finite subgroup $G\subset \Cr_n(\Bbbk)$
is conjugate to a biregular action on some rational projective variety $X$.
For any \emph{biregular} action of a finite cyclic 
group $G= \langle\tau \rangle$ on a smooth projective threefold $X$, one can define a 
subvariety $\F(\tau)\subset X$, a non-uniruled codimension one component of
the locus of fixed points.
The birational type of $\F(\tau)$ (if it is non-empty) does not depend on the choice 
of our model $X$, so it is an invariant of the subgroup $G\subset \Cr_n(\Bbbk)$
up to conjugacy. 

For $n=2$, the invariant $\F(\tau)$ is sufficient to distinguish
conjugacy classes of involutions $\tau$ in $\Cr_2(\Bbbk)$:

\begin{mtheorem}{\bf Theorem \cite{Bayle-Beauville-2000}.} \label{Theorem-involutions-Cr2}
Every element $\tau\in \Cr_2(\Bbbk)$ of order $2$ is conjugate to one
and only one of the following:
\begin{enumerate}
\item 
a linear involution acting on $\PP^2$, $\F(\tau)=\emptyset$;
\item 
\emph{de Jonqui\`eres} involution, $\F(\tau)$ is a hyperelliptic \textup(or elliptic\textup) curve of genus $g\ge 1$;
\item 
a \emph{Geiser} involution, $\F(\tau)$ is a non-hyperelliptic curve of genus $3$;
\item 
a \emph{Bertini} involution, $\F(\tau)$ is a 
non-hyperelliptic curve of genus $4$ whose canonical model lies on a singular quadric.
\end{enumerate}
Moreover, two elements $\tau,\, \tau'\in \Cr_2(\Bbbk)$ of order $2$
are conjugate if and only if $\F(\tau)\simeq \F(\tau')$.
In particular, $\tau\in \Cr_2(\Bbbk)$ is linearizable
\textup(i.e. conjugate to a linear involution acting on the projective plane\textup)
if and only if $\F(\tau)=\emptyset$.
\end{mtheorem}

In higher dimensions the last assertion of the above theorem 
is not true:
$\F(\tau)$ does not distinguish the conjugacy classes 
(see Example \ref{Example-non-rational-quotient}).
For example the field of invariants $\Bbbk(x_1,x_2,x_3)^\tau$ need not be rational.
Thus the isomorphism class of this field is another invariant 
distinguishing conjugacy classes of involutions.

In this paper we give a rough classification of elements of order two 
with $\F(\tau)\neq \emptyset$
in groups of 
birational selfmaps of not only rational  but  arbitrary
rationally connected varieties (not necessarily rational).
Our main result is the following.

\begin{mtheorem}{\bf Theorem.} 
\label{Theorem-main}
Let $Y$ be a three-dimensional rationally connected variety
and let $\tau\in \Bir(Y)$ be an element of order $2$ such that
$\F(\tau)\neq \emptyset$.
Then $\tau$ is conjugate to one of the following actions on a threefold $X$
\textup(birational to $Y$\textup):
\begin{enumerate}
\item[\type{C}]
$X$ is given by the equation
$\phi_0(u,v)x_0^2+\phi_1(u,v)x_1^2+\phi_2(u,v)x_3^2=0$ in 
$\PP^2_{x_0,x_1,x_2}\times \A^2_{u,v}$ and $\tau$ acts via 
$(x_0,x_1,x_2,u,v)\mapsto (-x_0,x_1,x_2,u,v)$.
\item[\type{D}]
$X$ is a $\tau$-equivariant del Pezzo fibration over a rational curve 
with trivial action of $\tau$.
The action on the generic fiber $X_\eta$ satisfies $\Pic(X_\eta)^{\tau}\simeq \ZZ$ and 
is one of the following:
\begin{enumerate}
\item
$K_{X_\eta}^2=1$ and $\tau_\eta: X_\eta\to X_\eta$ is the Bertini involution;
\item
$K_{X_\eta}^2=2$ and $\tau_\eta: X_\eta\to X_\eta$ is the Geiser involution;
\item
$1\le K_{X_\eta}^2\le 4$ and for the induced action 
$\tau_\eta: X_\eta\to X_\eta$ 
we have $\Fix(\tau_\eta,X_\eta)=C_\eta\cup \Lambda_\eta$, where 
$C_\eta$ is an elliptic curve, $C_\eta\in |-K_{X_\eta}|$,
and $\Lambda_\eta$ is a $0$-cycle of degree $4-K_{X_\eta}^2$;
\end{enumerate}

\item[\type{F^q}] 
$X$ is a Fano threefold with terminal $G\QQ$-factorial
singularities, $\Pic(X)^{\tau}\simeq \ZZ$, and one of the following holds
\begin{enumerate}
\item
$\dim |-K_X|\le 0$;
\item
$\dim |-K_X|=1$ and 
$|-K_X|$ defines a birational structure of a K3 fibration 
over $\PP^1$ \textup(see \xref{Proposition-Q-Fano-crepant}\textup);
\item
$\dim |-K_X|= 2$ and 
$|-K_X|$ defines a birational structure of an elliptic curve fibration 
over $\PP^2$ \textup(see \xref{Proposition-Q-Fano-crepant}\textup).
\end{enumerate}

\item[\type{F^c}] 
$X$ is a Fano threefold with 
canonical Gorenstein singularities:
\begin{enumerate}
\item
$X=X_6\subset \PP(1,1,1,1,3)$ is given by the equation
$y^2=\phi_6(x_1,\dots,x_4)$, the singularities of $X$
are cDV,  and $\tau$ acts via
\begin{enumerate}
\item
$y\mapsto -y$, or
\item
$x_1 \mapsto -x_1$ \textup(and $x_1$ appears in $\phi_6$ in even degrees only\textup).
\end{enumerate}
\item
$X$ is a double cover $X$ of a smooth quadric $W_2\subset \PP^4$ 
branched over a surface $B\subset W_2$ of degree $8$, the singularities of $X$
are terminal and $\tau$ is either
\begin{enumerate}
\item
the Galois involution of $X\to W_2$ or
\item
as in \xref{example-hyperelliptic-double-quadric-linear};
\end{enumerate}
\item
$X\subset \PP^4$ is a quartic with terminal singularities, $\tau$ is as in 
\xref{example-hyperelliptic-quartic};
\item
$X\subset \PP^5$ is a smooth intersection of a quadric and a cubic cone, $\tau$ is as in \xref{example-hyperelliptic-V6};
\item
$X\subset \PP^6$ is a smooth intersection of three quadrics, $\tau$ is as in \xref{example-hyperelliptic-V8}.
\end{enumerate}\end{enumerate}
\end{mtheorem}

\begin{case}{\bf Remark.} 
In all cases except for \type{F^q} we can describe the surface $\F(\tau)$:
\begin{center}
\begin{tabular}{l|p{287pt}} 
\multicolumn{1}{c|}{type} &\multicolumn{1}{c}{$\F(\tau)$ (up to birational equivalence)}
\\[4pt]
\hline
\type{C} & generically double cover of $\PP^2$
\\
\type{D}(a)&
fibration over $\PP^1$
into non-hyperelliptic curves of genus $4$ whose canonical model is contained in a singular quadric
\\
\type{D}(b) &
fibration over $\PP^1$
into non-hyperelliptic curves of genus $3$
\\
\type{D}(c) &
elliptic curve  fibration over $\PP^1$
\\
\type{F^c}(a)(i) & $\{\phi_6=0\}\subset \PP^3$
  \\
\type{F^c}(b)(i) & $B=W_2\cap U_4\subset \PP^4$, where $U_4$ is a quartic
  \\
\type{F^c} other cases&  K3 surface
\end{tabular}
\end{center}
\end{case}

\begin{case}{\bf Remark.}
Examples given in the corresponding sections show that all the
cases  \type{C}, \type{D}(a)-(c), \type{F^q}, (a)-(c), 
\type{F^c}(a)(i)-(e) really occur.
\end{case}

Note that our classification is really ``rough''.
First of all we do not provide detailed description, especially,
in the case \type{F^q}. The reason is the lack of a classification of 
singular Fano threefolds.
In fact, the involutions of type \type{F^q}(a) seems to be most difficult
to investigate.
Next, our cases can overlap and at the moment we cannot control this.
Finally, in many cases we cannot select \textit{rational} varieties from our list
to deduce a classification of involutions in the \textit{Cremona group}.
For example,  a general member $X$ as in \type{F^c} (a)-(d) and 
any member as in \type{F^c}(e)
is not rational \cite{Beauville1977},
\cite{Iskovskikh1980}. However there are a lot of examples of specific (singular) 
 Fano threefolds as in \type{F^c} (a)-(c) which are rational.

The paper is organized as follows.
Section \ref{section-Preliminaries} is preliminary.
In Section \ref{section-The-locus-of-fixed-points} we prove 
some easy technical facts concerning the set $\F(\tau)$.
In Sections \ref{section-Conic-bundles}, \ref{section-Del_Pezzo}, \ref{section-Non-Gorenstein},
and \ref{section-Gorenstein}
we prove Theorem \ref{Theorem-main} in cases \type{C}, \type{D}, \type{F^q}, and \type{F^c},
respectively. Section \ref{section-Gorenstein-examples} provides some examples 
of involutions of type \type{F^c}.

\section{Preliminaries}
\label{section-Preliminaries}
\begin{case}{\textbf{Notation.}}
\label{Notation}
\begin{itemize}
\item[] $\PP(a_1,\dots, a_n)$ denotes the weighted projective space,
\item[] $X_d\subset \PP(a_1,\dots, a_n)$ denotes a hypersurface of weighted degree $d$;
\item[] $X_{d_1\cdot d_2 \cdots d_r}\subset \PP(a_1,\dots, a_n)$ denotes a weighted complete intersection
of type $(d_1, d_2, \cdots, d_r)$.
\end{itemize}
We work over an algebraically closed field $\Bbbk$ of characteristic $0$.
When we say that a variety has, say, terminal singularities it means that
the singularities are not worse than that.
\end{case}

\begin{case} {\bf $G$-varieties.}\label{definition-GQ-Fano}
Let $G$ be a finite group.
In this paper a \emph{$G$-variety} is an algebraic variety $X$ provided with 
a biregular action of a finite group $G$.
Any morphism (resp. rational map) between 
$G$-varieties is usually supposed to be $G$-equivariant.
We say that a normal $G$-variety 
$X$ is \emph{$G\QQ$-factorial} if any $G$-invariant Weil 
divisor on $X$ is $\QQ$-Cartier.

A \textit{$G$-Fano-Mori fibration} (or \textit{$G$-Fano-Mori fiber space}) 
is a projective $G$-equivariant morphism $f: X\to Z$ such that
$f_*\OOO_X=\OOO_Z$, $\dim Z<\dim X$, $X$ has only terminal $G\QQ$-factorial
singularities, the anti-canonical 
divisor $-K_{ X}$ is ample over $Z$, and the relative $G$-invariant Picard 
number $\rho( X/Z)^G$ is one.
In the case $\dim X=3$, we have the following possibilities:
\begin{enumerate}
\item[\type{C}] 
$Z$ is a rational surface and the generic fiber $X_\eta$ is a conic;
\item[\type{D}] 
$Z\simeq \PP^1$ and the generic fiber $X_\eta$ is a smooth del Pezzo surface;
\item[\type{F}] 
$Z$ is a point and $ X$ is a, so-called, \textit{$G\QQ$-Fano threefold}.
\end{enumerate}
In these situations we say that $X/Z$ is of type 
\type{C}, \type{D}, \type{F}, respectively.
\begin{scase} {\bf Remark.}
\begin{enumerate}
\item 
In the case \type{C} there exists a non-empty Zariski open subset 
$U\subset Z$ such that the restriction $f_U: X_U\to U$ 
is a conic bundle. Similarly, in the case \type{D} 
there exists 
$U\subset Z$ such that $f_U: X_U\to U$ 
is a smooth del Pezzo fibration. 
\item 
We also will consider a class of Fano $G$-varieties different from \type{F}:
\begin{itemize}
 \item[\type{F^c}] 
Fano threefolds with canonical Gorenstein (not necessarily $G\QQ$-factorial) singularities.
\end{itemize}
In some situations this is more convenient for classification purposes.
Subclass of \type{F} consisting of $G\QQ$-Fano threefolds $X$
such that $K_X$ is not Cartier we denote by \type{F^q}.
\end{enumerate}
\end{scase}
\end{case}

\begin{mtheorem}{\bf Proposition.}\label{main-reduction}
Let $G$ be a finite group and let 
$X$ be a rationally connected $G$-variety.
Then there exists a $G$-Fano-Mori fibration $f: \bar X\to Z$ 
and a $G$-equivariant birational map $X\dashrightarrow \bar X$.
\end{mtheorem}
\begin{proof}[Outline of the proof]
Standard arguments (see e.g. \cite{Prokhorov2009e}) show that 
we can replace $X$ with a non-singular projective model so that
the action of $G$ on $X$ is biregular.
Next, we run 
the $G$-equivariant minimal model program ($G$-MMP): $X\dashrightarrow \bar X$ (in higher dimensions 
we can run $G$-MMP \emph{with scaling} \cite[Corollary 1.3.3]{BCHM}).
Note that \cite{BCHM} deals with varieties without group actions but 
adding an action of a finite group does not make a big difference 
(see \cite[Example 2.18]{Kollar-Mori-19988}).
Running this program we stay in the category of projective varieties 
with only terminal $G\QQ$-factorial singularities. 
Since $X$ is rationally connected, on the final step
the canonical divisor  cannot be nef \cite[Theorem 1]{Miyaoka1986} and so we get 
a $G$-Fano-Mori fibration $f: \bar X\to Z$. 
\end{proof}
\begin{scase}{\bf Definition.}
Notation as in \xref{main-reduction}.
If  $\dim X=3$, 
we say that our original $G$-variety belongs to class
\type{C} (resp. \type{D}, \type{F}) if 
so does $\bar X/Z$. Clearly, the choice of the birational map $X\dashrightarrow\bar X$  is not unique,
so the classes \type{C}, \type{D}, \type{F} can overlap. 
\end{scase}

\begin{case}
\label{pairs}
\textbf{$G$-minimal model program for pairs \cite{Alexeev-1994ge}.}
Let $X$ be a normal $G$-variety and let $\MMM$ be an invariant linear system
of Weil divisors without fixed components. 
In this situation we say that $(X,\MMM)$ is a \textit{$G$-pair}.
In the standard way one can define the discrepancy 
$a(E,X,\MMM)$ of a prime divisor $E$ with respect to $(X,\MMM)$.
We say that $(X,\MMM)$ is \textit{terminal} (resp. \textit{canonical}) if 
$a(E,X,\MMM)>0$ (resp. $\ge 0$) for all exceptional divisors $E$ over $X$.
One can run $G$-minimal model program in the  category of $G$-pairs 
$(X,\MMM)$ such that 
$X$ is $G\QQ$-factorial and $(X,\MMM)$ is terminal.
In particular, for any $G$-pair $(X,\MMM)$ there exists a \textit{terminal model }
$f: (X',\MMM')\to (X,\MMM)$, where $f$ is a $G$-equivariant morphism and 
$(X',\MMM')$ is a terminal $G$-pair such that 
$K_{X'}+\MMM'$ is $f$-nef,
$\MMM'$ is the birational transform of $\MMM$, and $X'$ is $G\QQ$-factorial.
Moreover, we can write 
\[
K_{X'}+\MMM'=f^*(K_X+\MMM) - \sum a_i E_i, 
\]
where $E_i$ are $f$-exceptional divisors, $a_i\ge 0$ for all $i$,
and $a_i= 0$ for all $i$ if and only if $(X,\MMM)$ is canonical.
In the last case we say that $f$ is \textit{log crepant}.
\end{case}

\begin{case}{\bf Varieties of minimal degree.}
For convenience of references we recall the following well-known fact.
\begin{stheorem}{\bf Theorem.}
\label{Theorem-Del-Pezzo-Bertini}
Let $W\subset \PP^N$ be a $n$-dimensional 
projective variety not lying in a hyperplane.
Then $\deg W\ge \codim W +1$ and the equality holds if and only if 
$W$ is one of the following:
\begin{enumerate} 
\renewcommand\labelenumi{(\arabic{enumi})}
\renewcommand\theenumi{(\arabic{enumi})}
 \item 
 $W=\PP^N$;
 \item
$W=W_2\subset \PP^N$ is a quadric; 
 \item
$W$ is the image of $\PP_{\PP^2}(\EEE)$, where 
$\EEE=\OOO_{\PP^2}(2)\oplus \bigoplus_{i=1}^{n-2} 
\OOO_{\PP^2}$ under the \textup(birational\textup) morphism given by $|\OOO(1)|$;
 \item
$W$ is the image of $\PP_{\PP^1}(\EEE)$, where 
$\EEE=\bigoplus_{i=1}^{n} 
\OOO_{\PP^1}(a_i)$, $a_i\ge 0$ under the \textup(birational\textup) morphism given by $|\OOO(1)|$.
\end{enumerate}
\end{stheorem}
\end{case}

\begin{case}{\bf Fano threefolds.}\label{Fano-threefolds}
Let $X$ be a Fano threefold with only canonical Gorenstein singularities.
By the Riemann-Roch formula and Kawamata-Viehweg vanishing 
$\dim |-K_X|=g+1$, where $g=\g(X)$ is a positive integer such that $-K_X^3=2g-2$.
This number is called the \textit{genus} of $X$. 
Let $\Phi: X \dashrightarrow \PP^{g+1}$ be the anti-canonical map.

All such Fanos are divided in the following groups
(see \cite{Iskovskikh-1980-Anticanonical}, \cite{Iskovskikh-Prokhorov-1999}, \cite{Jahnke-Radloff-2006},
\cite{Przhiyalkovskij-Chel'tsov-Shramov-2005en}):
\begin{itemize}
\item 
$\Bs |-K_X|\neq \emptyset$ and $\Phi(X)$ is a surface as in \xref{Theorem-Del-Pezzo-Bertini};
\item 
(hyperelliptic case)
$\Bs |-K_X|= \emptyset$, $\Phi$ is a double cover of a threefold $W\subset \PP^{g+1}$
as in \xref{Theorem-Del-Pezzo-Bertini};
\item 
(trigonal case)
$\Phi$ is an embedding and the image $\Phi(X)$ is not an intersection of quadrics,
in this case quadrics passing through $\Phi(X)$ cut out a fourfold as in \xref{Theorem-Del-Pezzo-Bertini};
\item
(main series)
$\Phi$ is an embedding and the image $\Phi(X)$ is an intersection of quadrics.
\end{itemize}
In the above notation,  $X$ is called 
\textit{del Pezzo threefold} if its anti-canonical class is divisible by $2$ in 
$\Pic(X)$, i.e. $-K_X\sim 2A$ for some $A\in \Pic(X)$. 
The \textit{degree} of $X$ is defined as $d=\dd(X):=A^3$.
\end{case}

\begin{mtheorem}{\bf Theorem 
(\cite[Ch. 2, \S 1]{Iskovskikh-1980-Anticanonical}, \cite[Ch. 1, \S 6]{Fujita-all}).}
\label{th-del-Pezzo-2}
Let $X$ be a del Pezzo threefold and let $A:=-\frac{1}{2}K_X$.
Then $\dim |S|=\dd(X)+1$ and $\dd(X)\le 8$. Moreover,
\begin{enumerate}
\item
If $d=1$, then $X\simeq X_6\subset \PP(1^3,2,3)$. The linear system $|2A|$ 
defines a double cover $X\to \PP(1^3,2)$ whose 
branch locus $B\subset \PP(1^3,2)$ is a surface of weighted degree $6$
such that the pair $(\PP(1^3,2),\frac 12 B)$ is klt.

\item
If $d=2$, then $X\simeq X_4\subset \PP(1^{4},2)$. 
The linear system $|A|$ defines a double cover $X\to \PP^3$ whose 
branch locus $B\subset \PP^3$ is a surface of degree $4$ 
such that the pair $(\PP^3,\frac 12 B)$ is klt.

\item
If $d=3$, then $X\simeq X_3\subset \PP^{4}$.
\item
If $d=4$, then $X\simeq X_{2\cdot 2}\subset \PP^{5}$.
\end{enumerate}
\end{mtheorem}

\section{The locus of fixed points}
\label{section-The-locus-of-fixed-points}
\begin{case}\label{Notation-The-locus-of-fixed-points}
From now on  $G$ denotes a finite cyclic group generated by 
an element $\tau$.
Let $X$ be a $G$-variety of dimension $n$.
Denote the set of fixed points by $\Fix(\tau,X)$. 
We will write $\Fix(\tau)$ instead 
if no confusion is likely.
\end{case}

\begin{case}{\bf Remark.}\label{Remark-1}
Assume that $X$ has only $G\QQ$-factorial terminal singularities
and let $F=\cup F_i$ be the union of all the $(n-1)$-dimensional components $F_i\subset \Fix(\tau)$.
Then $\Sing(F)\subset \Sing(X)$, $F$ is is a disjoint union of its irreducible components and
each component $F_i$ is normal. 
Indeed, for every point $P\in F$, the induced action of $\tau$ 
on the Zariski tangent space $T_{P,F}$ is trivial.
On the other hand, the action is non-trivial on $T_{P,X}$. 
Hence $T_{P,F} \subsetneq T_{P,X}$ and so $\Sing(F)\subset \Sing(X)$.
In particular, distinct components $F_i, \, F_j\subset \Fix(\tau)$
cannot meet each other outside of $\Sing(X)$.
Note that $\codim \Sing(X)\ge 3$ because the singularities of $X$ are terminal.
On the other hand, $F_i$ and $F_j$ are $\QQ$-Cartier divisors
because $X$ is $G\QQ$-factorial.
Therefore, $F_i\cap F_j=\emptyset$.
The same arguments show that each $F_i$ is 
smooth in codimension one. Since $F_i$ is Cohen-Macaulay 
(see e.g. \cite[Corollary 5.25]{Kollar-Mori-19988}), it is normal.
\end{case}

\begin{case}{\bf Remark.}\label{Remark-2}
 Assume that $X$ has a structure of a $G$-Mori-Fano fiber space $f: X\to Z$
and assume that there exists a non-uniruled irreducible component $S\subset \Fix(\tau,X)$
of dimension $n-1$. Then $S$ dominates $Z$ (see e.g. \cite[Corollary 1.3]{Hacon2007}).
Therefore, the action of $\tau$ on $Z$ is trivial.
\end{case}

\begin{mtheorem} {\bf Proposition.}\label{Proposition-components}
Assume that $X$ is rationally connected.
\begin{enumerate}
 \item 
There exists at most one irreducible component $S\subset \Fix(\tau,X)$
of dimension $n-1$ which is not uniruled.
 \item 
Let $\psi: X\dashrightarrow X'$ be a birational map where $X'$ 
is a projective variety with Kawamata log terminal singularities
\textup(with respect to some boundary\textup).
Assume that there is a divisorial component $S\subset \Fix(\tau,X)$ which is not uniruled.
Then $\psi_* S$ is a divisor on $X'$.
\end{enumerate}
\end{mtheorem}

\begin{proof}
(i) Assume the contrary: there are two such components $S_1$ and $S_2$.
As in  \ref{main-reduction} we may assume that $X$ is projective and smooth.
Run the equivariant MMP $X \dashrightarrow \bar X$. 
Non-uniruled divisors $S_i$ cannot be contracted (see e.g. \cite[Corollary 1.3]{Hacon2007}).
Therefore, their images $\bar S_1$ and $\bar S_2$ are also divisors contained in $\Fix(\tau,\bar X)$.
Since $X$ is rationally connected, it cannot have a minimal model, i.e. 
the canonical divisor $K_{\bar X}$ cannot be nef \cite[Theorem 1]{Miyaoka1986}.
Thus $\bar X$ has a structure of a $G$-Fano-Mori fiber space $f: \bar X\to Z$.
By Remark \ref{Remark-1} the divisors $\bar S_1$, \, $\bar S_2$ 
do not meet each other. By Remark \ref{Remark-2}
both $\bar S_1$ and $\bar S_2$ are $f$-ample and the action of $G$ on $Z$ is trivial. 
Restricting $\bar S_1$ and $\bar S_2$ to a general fiber $F=f^{-1}(z)$ we 
get two ample disjointed divisors. Hence, $\dim F=1$ and so  
$F\simeq \PP^1$.
Since the action is non-trivial on $F$, we see that
$\Fix(\tau,\bar X)\cap F$ consists of two points $\bar S_1\cap F$ and $\bar S_2\cap F$.
On the other hand, $Z$ is rationally connected and so 
$\bar S_i$ cannot be a section of $f$ generically,
a contradiction.

(ii) follows by \cite[Cor. 1.6]{Hacon2007}.
\end{proof}

\begin{case}{\bf Notation.}
We denote the non-uniruled codimension one 
component of $\Fix(\tau,X)$ by $\F(\tau,X)$ or simply by $\F(\tau)$ (if 
no confusion is likely).
By Proposition \ref{Proposition-components} (ii) the birational type of 
$\F(\tau,X)$ does not depend on the choice of the projective model $X$
(with only log terminal singularities)
on which $G$ acts biregularly. In particular, we can define the Kodaira dimension
$\kappa(\tau,X)$ to be equal to $\kappa(\F(\tau,X))$ if $\F(\tau,X)\neq \emptyset$
and $\kappa(\tau,X)=-\infty$ otherwise. 
\end{case}

\section{Conic bundles}
\label{section-Conic-bundles}
\begin{case}{\bf Notation.}
From now on, let $G=\{1,\tau\}$ be a group of order $2$.
\end{case}

\begin{mtheorem} {\bf Lemma.}\label{Lemma-Conic-bundles}
Let $f: X\to Z$ be a $G$-Fano-Mori fiber space over 
a rational surface $Z$. Assume that the action of $G$ 
on $Z$ is trivial. Then $f$ is $G$-birationally equivalent to
the action on the hypersurface in $\PP^2_{x_0,x_1,x_2}\times \bA^2_{u,v}$ given by 
\begin{equation}\label{equation-1-Conic-bundles}
\phi_{0}(u,v)x_0^2+ \phi_{1}(u,v)x_1^2+\phi_{2}(u,v)x_2^2=0
\end{equation}
via
\begin{equation}\label{equation-2-Conic-bundles}
\tau: (x_0,x_1,x_2; u,v) \longmapsto (-x_0,x_1,x_2; u,v).
\end{equation}
Moreover, the quotient $X/G$ is rational.
\end{mtheorem}

\begin{proof}
First, by shrinking $Z$ we may assume that $Z\subset \bA^2_{u,v}$
is a Zariski open subset, $X$ is smooth, and $f$ is a conic bundle over $Z$.
Then the linear system $|-K_X|$ defines a $G$-equivariant embedding 
$X \hookrightarrow \PP^2_{x_0,x_1,x_2}\times \bA^2_{u,v}$ so that the fibers of 
$f$ are conics in $\PP^2_{x_0,x_1,x_2}$ (see \cite[Proposition 1.2]{Beauville1977}). 
Clearly, the action on $\PP^2\times \bA^2$ can be written in the form 
\eqref{equation-2-Conic-bundles}.
Then by changing the coordinate system we may assume that $X$ is given by 
\eqref{equation-1-Conic-bundles}.
Finally, $X/G$ is rational because 
$\Bbbk(X/G)\simeq \Bbbk(u,v,x_1/x_2)$.
\end{proof}

\begin{stheorem}{\bf Corollary.}
Let $Y$ be a rationally connected 
$G$-threefold with $G=\{1,\tau\}$ and $\kappa(\tau,Y)\ge 0$. 
Assume that $(Y,G)$ belongs to the class \type{C}
and let $f: X\to Z$ be a conic bundle model, that is, 
a $G$-Fano-Mori fibration with $\dim Z=2$ such that there exists a 
$G$-equivariant birational map $Y\dashrightarrow X$.
Then the action on $Z$ is trivial and the 
restriction $f_{\F(\tau)}: \F(\tau)\to Z$
is a generically double cover.
Therefore, the action of $G$ is conjugate to that 
given by \eqref{equation-1-Conic-bundles}-\eqref{equation-2-Conic-bundles}.
\end{stheorem}
\begin{proof}
Let $S:=\F(\tau)$. By Remark \ref{Remark-2} the action on $Z$ is trivial and $f(S)=Z$.
This proves the first statement.
Further, for a general fiber $F$ we have $F\simeq \PP^1$ and
the action of $G$ on $F$ is not trivial. 
Then $F\cap S$ consists of two points and so $f|_S: S\to Z$ is a generically double cover.
\end{proof}

\begin{case} {\bf Remark.}
Recall that the \textit{irrationality degree} of an algebraic variety $V$
is the minimal degree
of a dominant rational map from $V$ to $\PP^n$ \cite{Yoshihara-1994}.
By the above corollary the irrationality degree of $\F(\tau)$,
where $G$ is an involution with $\F(\tau)\neq \emptyset$ acting on a conic bundle,
equals to $2$. The following construction shows that
\textit{any} algebraic surface of irrationality degree $2$ is birational to
$\F(\tau)$ for some involution in $\Cr_3(\Bbbk)$.

\begin{scase} {\bf Construction.}
Let $\pi:S \dashrightarrow Z$ be a generically finite rational map of degree $2$,
where $Z$ is a rational surface.
By shrinking $Z$ we may assume that $\pi$ is a morphism, $Z$ is a Zariski open subset in 
$\bA^2_{u,v}$, and $S$ is given in $\bA^2_{u,v}\times \PP^1_{x}$ by the equation
\[
x^2=\phi(u,v). 
\]
Now define the action of $\tau$ on 
$X:=\bA^2_{u,v}\times \PP^1_{x}$ as follows:
\[
\tau : (u,\, v,\, x)\longmapsto \Bigl(u,\, v,\, \phi(u,v)x^{-1}\Bigr).
\]
Clearly, $\Fix(\tau)=S$.
\end{scase}
\end{case}

\begin{case}
Finally, we give a series of examples of involutions
in groups of birational selfmaps of conic bundles with $\F(\tau)= \emptyset$.
In particular, this gives examples of non-linearizable involutions in 
$\Cr_3(\Bbbk)$ with $\F(\tau)= \emptyset$ (cf. Theorem \ref{Theorem-involutions-Cr2}).

\begin{scase} {\bf Construction.}
Let $h:Y\to Z$ be a conic bundle over a surface.
Assume that there exists 
a rational surface $S\subset Y$ such that 
the restriction $h_S: S\to Z$ is a generically finite morphism of degree $2$
(so $S$ is a $2$-section of $h$ generically).
By shrinking $Z$ we may assume that $h_S$ is, in fact, finite. 
Consider the base change $f: X=Y\times_Z S\to S$. Then $f$ is a locally trivial
rational curve fibration. In particular, $X$ is rational.
Then the Galois involution $\tau$ of the double cover $X\to Y$ defines an element of 
$\Cr_3(\Bbbk)$. 
Note that if $Y=X/\tau$ is not rational, then $\tau$ is not linearizable.
Indeed, otherwise $X/\tau$ is birationally equivalent to the rational variety $\PP^3/\tau'$,
where $\tau'$ is a linear involution of $\PP^3$ that conjugate to $\tau$.
\end{scase}

\begin{scase} {\bf Example.}\label{Example-non-rational-quotient}
Let $V\subset \PP^4$ be a smooth cubic hypersurface
and let $L\subset V$ be a line.
Let $Y\to V$ be the blowup of $L$ and let $S$ be the exceptional divisor.
The projection away from $L$ induces a conic bundle structure $Y\to \PP^2$
and $S$ is a rational $2$-section.
The cubic $V$ is not rational \cite{Clemens-Griffiths}.
According to Lemma \ref{Lemma-Conic-bundles}
our construction gives an example of an involution in $\Cr_3(\Bbbk)$
with $\F(\tau)=\emptyset$ which is not linearizable.
\end{scase}
\end{case}

\section{Del Pezzo fibrations}
\label{section-Del_Pezzo}
In this section we study the case \type{D}. 

\begin{mtheorem}{\bf Lemma.}
Let $\Gamma$ be any finite group.
There is a natural bijective correspondence between 
\begin{itemize}
\item  
the set of all $\Gamma$-del Pezzo fibrations $f: X\to \PP^1$ with trivial action on $\PP^1$
modulo $\Gamma$-equivariant birational transformations over $\PP^1$, and
\item 
the set of all del Pezzo surfaces $V$  over $\Bbbk(t)$ 
admitting an effective $\Gamma$-action  with
$\Pic(V)^\Gamma\simeq\ZZ$ modulo $\Gamma$-isomorphisms.
\end{itemize}
\end{mtheorem}

\begin{proof}[Outline of the proof]
Let $\Bbbk(\PP^1)=\Bbbk(t)$ and let $\eta=\Spec \Bbbk(t)$ be the 
generic point.
Let $f: X\to \PP^1$ be a $\Gamma$-del Pezzo fibration with   trivial action on $\PP^1$.
Then $\Gamma$ acts on  the generic fiber $V:=X_\eta$ 
which is a del Pezzo surface over $\Bbbk(t)$ with $\Pic(V)^\Gamma\simeq\ZZ$.
Conversely, if $\Gamma$ acts on a del Pezzo surface $V$ defined 
over $\Bbbk(t)$ so that $\Pic(V)^\tau\simeq\ZZ$, we can take a subring 
$\Bbbk\subset R\subset \Bbbk(t)$ such that the quotient 
field of $R$ coincides with $\Bbbk(t)$, $Z:=\Spec R$ is an rational affine 
normal  curve, and
there exists a variety $X_Z$ over $R$ such that
$X_Z\times _{\Spec R} \Spec \Bbbk(t)=V$. 
Then as in \ref{main-reduction} 
we get a $\Gamma$-del Pezzo fibration $X/\PP^1$.
\end{proof}
Recall that $G=\{1,\, \tau\}$ in our case.
\begin{mtheorem} {\bf Proposition.}\label{Theorem-del-Pezzo}
Let $f:X\to Z$ be a $G$-equivariant 
del Pezzo fibration with $\Pic(X/Z)^G\simeq \ZZ$ and $\F(\tau)\neq \emptyset$
over a rational curve $Z$.
Let $X_\eta$ be the generic fiber and let $d:=K_{X_\eta}^2$.
Then one of the following holds:
\begin{enumerate}
\item \label{Theorem-del-Pezzo-Bertini}
$d=1$, $\tau$ is a Bertini involution on $X_\eta$;
\item\label{Theorem-del-Pezzo-Geiser}
$d=2$, $\tau$ is a Geiser involution on $X_\eta$;
\item
$1\le d\le 4$ and $\tau$ fixes pointwisely an elliptic curve on $X_\eta$.
\end{enumerate}
Moreover, $\Pic(X_\eta/G)\simeq\ZZ$.
\end{mtheorem}

\begin{proof}
Recall that  by Remark \ref{Remark-2} the 
action on $Z$ is trivial.
Let $S:=\F(\tau,X)$. 
Let $X_\eta$ be the generic fiber and let $S_\eta\subset X_\eta$
be the generic fiber of $S/Z$.
Thus $S_\eta$ is a curve of $\tau$-fixed points. 
Since $\kappa(S)\ge 0$, $p_a(S_\eta)>0$.
If $p_a(S_\eta)>1$, then by \cite{Bayle-Beauville-2000} (cf. Theorem \ref{Theorem-involutions-Cr2})
we have cases (i) and (ii). Assume that $p_a(S_\eta)=1$. Then
$S_\eta$ is an elliptic curve and $S_\eta\in |-K_{X_\eta}|$.
In this case the statement follows by Lemma \ref{Lemma-del-Pezzo-c1}
below. The last assertion is clear because $\rho(X_\eta/G)=\rho(X_\eta)^G=1$.
\end{proof}

\begin{stheorem}{\bf Lemma.}\label{Lemma-del-Pezzo-c1}
Let $V$ be a del Pezzo surface over a field $\LL$ of characteristic $0$, 
let $\tau\in \Aut_{\LL}(V)$ be an element of order $2$, and let
$W:=V/\tau$.
Assume that $\F(\tau, V)$ is a \textup(smooth\textup) elliptic curve. Then 
$W$ is a del Pezzo surface with only Du Val singularities of type $A_1$,
$K_W^2=2K_V^2$, and 
\[
\# \Sing (W)+ K_V^2=4. 
\]
In particular, $K_V^2\le 4$. 
\end{stheorem}
\begin{proof}
We may assume that $\LL=\CC$.
Let $\pi: V\to W$ be the quotient map. 
Put $s:=\# \Sing (W)$ and $d:=K_V^2$.
Near each fixed point $P\in V$ 
the action of $\tau$ in suitable (analytic) coordinates can be written as 
$\tau: (x_1,x_2) \mapsto (-x_1,x_2)$ or $(-x_1,-x_2)$.
In the first case $\pi(P)$ is smooth and in the second one $\pi(P)$ is of type $A_1$.
Let $C$ be the sum of one-dimensional components of $\Fix(\tau)$
and let $C_1:=\F(\tau, V)\subset C$. 

We claim that $C_1=C\sim -K_V$. Assume that  $C_1\not\sim -K_V$. Since 
$(K_V+C_1)\cdot C_1=2p_a(C_1)-2=0$ and $K_V\cdot C_1<0$, we have $C_1^2>0$. 
By the Hodge index theorem $(K_V+C_1)\cdot K_V=(K_V+C_1)^2<0$, so $d=K_V^2<-K_V\cdot C_1$.
On the other hand, the action of $\tau$ on the linear system $|-r_d K_V|$ is non-trivial,
where $r_1=2$ and $r_d=1$ if $d>1$.
Indeed, otherwise the morphism given by $|-r_dK_V|$ passes through $W$ and then
$\tau$ must be either a Bertini or Geiser involution.
Hence $C$ is a component of a divisor $D\in |-r_dK_V|$.
Therefore, $d<-K_V\cdot C_1\le -K_V\cdot D = dr_d$. This is possible 
only if $d=1$ and $D=C_1$. But then $p_a(C_1)=2$, a contradiction.
Thus $C_1\sim -K_V$. 
Note that $C$ is smooth (because  $C\subset \Fix(\tau)$). 
In particular, $C$ is a disjoint union of its irreducible components.
Since $C_1\sim -K_V$ is ample, $C=C_1$.

Now by the Hurwitz formula $K_V=\pi^*K_W+C$, so
$2K_V=\pi^*K_W$. Hence $-K_W$ is ample, i.e. 
$W$ is a del Pezzo surface with only Du Val singularities of type $A_1$
and $K_W^2=2d$.
Applying the Noether formula to the minimal resolution of $W$ we get
\[
\rho(W)+s= 10-K_W^2=10-2d.
\]
On the other hand, by the Lefschetz fixed-point-formula we have
\[
s=2+ \operatorname{Tr}_{H^2(V,\mathbb R)}(\tau)=
2+\rho(W)- (\rho(V)-\rho(W))=
2+2\rho(W)- 10+d.
\]
Thus $s+d=4$.
\end{proof}
Now we give two examples of involutions in $\Cr_3(\Bbbk)$ 
of types \ref{Theorem-del-Pezzo}(i)-(ii).

\begin{case} {\bf Example (cf. \cite[7.7, 8.6, 4.2, 5.2]{Prokhorov2010a}).}
Let $X=X_4\subset \PP(1^4,2)_{x_1,\dots,x_4,y}$ be a hypersurface of degree $4$ 
having only terminal singularities and let $\pi: X \to \PP^3_{x_1,\dots,x_4}$ be the projection.
Then $\pi$ is a double cover branched over a quartic $B\subset \PP^3$
and $X$ is a del Pezzo threefold of degree $2$.
Let $\tau$ be the Galois involution. 
We call it the \emph{Geiser involution} of $X$. 
One can write the equation of $X$  in the form 
\[
y^2=\phi_4(x_1,\dots,x_4), 
\]
where $\phi_4=0$ is the equation of $B$ and $\tau$ acts via 
$x_i\mapsto x_i$, $y \mapsto -y$. 
Let $l_1$ and $l_2$ are general linear forms in $x_1,\dots,x_4$.
Then the map
\[
X \dashrightarrow \PP^1,\qquad
(x_1,\dots,x_4,y)\dashrightarrow (l_1,l_2)
\]
defines a birational $\tau$-equivariant structure 
of del Pezzo fibration on $X$ as in \ref{Theorem-del-Pezzo-Geiser} 
of Proposition \ref{Theorem-del-Pezzo}. 
Clearly, $\Fix(\tau)\simeq B$ and this surface is either K3
with only Du Val singularities or birationally ruled.
For a general choice of $B$, the variety $X$ is not rational.
However, for some special $B\subset \PP^3$ the variety $X$ can be rational.
For example, let $\tilde \PP^3\to \PP^3$ be the blowup of $6$ points in general position 
and let $\tilde \PP^3\to X$ be the contraction of $K$-trivial curves.
Then $X$ is a del Pezzo threefold as above. This $X$ is rational 
and has exactly $16$ singular points. Hence $B$
is a Kummer surface. See \cite{Prokhorov2010a} for more examples.
\end{case}

\begin{case} {\bf Example (cf. \cite[7.8,  5.2]{Prokhorov2010a}).}
As above, starting with a hypersurface 
$X=X_6\subset \PP(1^3,2,3)$ we get the \emph{Bertini involution} on $X$.
A general projection defines a birational $\tau$-equivariant structure 
of a del Pezzo fibration on $X$ as in \ref{Theorem-del-Pezzo-Bertini} 
of Proposition \ref{Theorem-del-Pezzo}. In this case $\Fix(\tau)$ 
is a hypersurface of degree $6$ in $\PP(1^3,2)$, a surface of general type.
\end{case}

\begin{case} \label{construction-del-Pezzo-}
Involutions 
of the type \ref{Theorem-del-Pezzo}(iii) can be given by the following
series of examples.
\end{case}

\begin{scase} {\bf Example.}\label{Example-del-Pezzo-C-4}
Let $Q\subset \PP^3_{\Bbbk(t)}$ be a smooth quadric defined over $\Bbbk(t)$.
Assume that $\Pic(Q)\simeq \ZZ$. For example we can take $Q=\{x_0x_1 +x_2^2+tx_3^2=0\}$.
Let $V\to Q$ be the double cover branched over a smooth curve $C\in |-K_Q|$
(also defined over $\Bbbk(t)$).
Then $V$ is a del Pezzo surface of degree $4$. 
It is easy to see that $V$ is an intersection of two quadrics $Q_1$, $Q_2$ in $\PP_{\Bbbk(t)}^4$
such that $Q_1$ is a cone over $Q$ and $Q_2\cap Q=C$.
Let $\tau$ be the Galois involution of $V/Q$.
Since $\Pic(Q)\simeq \ZZ$, we have $\Pic(V)^\tau\simeq \ZZ$.
We get an action of $\tau$ on a del Pezzo fibration as in (iii) of Proposition \ref{Theorem-del-Pezzo}.
\end{scase}

\begin{scase} {\bf Example.}\label{Example-del-Pezzo-C-3}
Let $V\subset \PP^3_{\Bbbk(t)}$ be a smooth cubic surface 
given by the equation $x_0^2x_1+ x_1\phi_2(x_1,x_2,x_3)+\phi_3(x_1,x_2,x_3)=0$,
where $\phi_d$ is a homogeneous polynomial of degree $d$ and $\phi_3(x_1,x_2,x_3)$
is irreducible over $\Bbbk(t)$. Thus $(1:0:0:0)$ is an Eckardt point
and lines passing through it are conjugate under $\Gal(\overline{\Bbbk(t)}/\Bbbk(t))$.
The involution $\tau$ acts via $x_0\mapsto -x_0$ and the quotient
$V/\tau$ is a del Pezzo surface of degree $6$ with a point of type $A_1$.
As above we get an action of $\tau$ on a del Pezzo fibration as in \ref{Theorem-del-Pezzo}(iii).
\end{scase}

 \begin{scase} {\bf Example.}\label{Example-del-Pezzo-C-2}
Let $V\subset \PP_{\Bbbk(t)}(1,1,1,2)$ 
be a smooth hypersurface of weighted degree $4$.
We can write its equation as $y^2=\phi_4(x_1,x_2,x_3)$.
Assume that $\phi_4$ contains terms of even degree in $x_1$ only.
Then $V$ is invariant under the involution $\tau: x_1\mapsto -x_1$.
The quotient $V/\tau$ is a hypersurface of weighted degree $4$
in $\PP_{\Bbbk(t)}(1,1,2,2)$. 
It is a del Pezzo surface of degree $4$ with two points of type $A_1$.
As above we get an action of $\tau$ on a del Pezzo fibration as in \ref{Theorem-del-Pezzo}(iii).
\end{scase}

\begin{scase} {\bf Example.}\label{Example-del-Pezzo-C-1}
Let $V\subset \PP_{\Bbbk(t)}(1,1,2,3)$ 
be a smooth hypersurface of weighted degree $6$.
We can write its equation as $z^2=\phi_6(x_1,x_2,y)$.
Assume that $\phi_6$ contains only terms of even degree in $x_1$.
The involution $\tau$ acts on $V$ via $\tau: x_1\mapsto -x_1$. 
As above we get an action of $\tau$ on a del Pezzo fibration
as in \ref{Theorem-del-Pezzo}(iii).
\end{scase}

\section{Non-Gorenstein Fano threefolds.}
\label{section-Non-Gorenstein}
In this section we study the case \type{F^q}.
Since we assert nothing if $\dim |-K_X|\le 0$, we 
may assume that $\dim |-K_X|\ge 1$. 
The essential result here is to show that for a Fano threefold with $\dim |-K_X|\ge 3$
admitting an involution with $\F(\tau)\neq \emptyset$ 
there is an equivariant birational transformation to 
a threefold with a structure of a fibration of type \type{C} or \type{D}, or
to a Fano threefold with \textit{canonical Gorenstein} singularities.
We need the following easy but very useful lemma.

\begin{mtheorem}{\bf Lemma.}\label{Lemma-index>1}
Let $X$ be a normal $G$-threefold and let $\MMM$ be a $G$-invariant 
linear system without fixed components
such that $-K_X-\MMM \qq D$, where $D$ is an effective $\QQ$-divisor.
If either $D\neq 0$ or the pair $(X,\MMM)$ is not canonical, 
then $(X,G)$ is of type \type{C} or \type{D}. 
\end{mtheorem}

\begin{proof}
Take a $G$-invariant pencil $\LLL\subset \MMM$ (without fixed components)
and replace $(X,\LLL)$ with its a terminal $G\QQ$-factorial 
model (see \ref{pairs}). Then we can write
$-K_X\qq \LLL+ D$, where $D>0$.
Run $G$-equivariant $(K_{X}+\LLL)$-MMP.
On each step we contract a $(K+ \LLL)$-negative extremal ray and $-(K+\LLL)\qq D$.
Hence $D$ is not contracted and  $K+ \LLL$ cannot be nef.
If $(X,G)$ is not of type \type{C} nor \type{D}, at the end we get 
a pair $(X',\LLL')$ such that $\Pic(X')^G\simeq\ZZ$ and $-K_{X'}\qq \LLL'+D'$,
where $D'>0$. 
We can write $-K_{X'}\qq a' \LLL'$ for some $a'>1$.
By our construction the pair $(X',\LLL')$ is terminal.
Hence a general member $L'\in \LLL'$ is smooth and contained in the smooth locus of 
$X'$ \cite[Lemma 1.22]{Alexeev-1994ge}. 
By the adjunction formula $-K_{L'}=(a'-1) L'|_{L'}$.
Hence $L'$ is a del Pezzo surface. In particular, $L'$ is rational
and so $X$ has a $G$-invariant  pencil rational surfaces.
Resolving the base locus and running the relative $G$-MMP
we get a $G$-Fano-Mori fibration of type \type{D} or \type{C}.
\end{proof}

\begin{stheorem}{\bf Corollary.}
 \label{Corollary-M1}
Let $X$ be a $G\QQ$-Fano threefold.
Assume that $-K_X\sim D+M$, where $D$ and $M$ are effective divisors such that 
$D\neq 0$ and $\dim |M|>0$. Then $(X,G)$ is of type \type{C} or \type{D}. 
\end{stheorem}

\begin{stheorem}{\bf Corollary (cf. \cite[\S 4]{Alexeev-1994ge}).}
\label{Corollary}
Let $X$ be a $G\QQ$-Fano threefold.
Assume that $\dim |-K_X|\ge 1$ and $(X,G)$ is not of type \type{C} nor \type{D}. 
Then the pair $(X,|-K_{X}|)$ is canonical. 
In particular, $|-K_X|$ has no fixed components.
\end{stheorem}

\begin{mtheorem}{\bf Proposition (cf. \cite[Theorem 4.5]{Alexeev-1994ge}).}
\label{Proposition-Q-Fano-crepant}
Let $X$ be a non-Gorenstein $G\QQ$-Fano threefold.
Assume that $\dim |-K_X|\ge 1$ 
and $(X,G)$ is not of type \type{C} nor \type{D}. 
Then there are $G$-equivariant birational morphisms 
\[
X \overset{f}\longleftarrow \tilde X \overset{h}\longrightarrow Y,
\]
where 
$f$ is a log crepant terminal model of the pair $(X,|-K_X|)$
and the singularities of $\tilde X$ are terminal Gorenstein.
Furthermore, there are the following possibilities:
\begin{enumerate}
\item 
$Y\simeq\PP^1$, $\dim |-K_X|=1$, and $h$ is a K3-fibration;
\item 
$Y\simeq\PP^2$, $\dim |-K_X|=2$, and $h$ is an elliptic fibration;
\item 
$Y$ is a Fano threefold with only Gorenstein canonical singularities and
$h$ is a birational $K$-crepant morphism.
\end{enumerate}
\end{mtheorem}

\begin{proof}[Proof \textup(cf. \textup{\cite[Theorem 4.5]{Alexeev-1994ge}}\textup).]
By Corollary \ref{Corollary}
the pair $(X,|-K_X|)$ is canonical.
Let $f$ be a log crepant terminal model of $(X,|-K_X|)$.
Thus $(\tilde X,|-K_{\tilde X}|)$ is terminal and $\tilde X$ is $G\QQ$-factorial.
Then $|-K_{\tilde X}|$ is base point free \cite[Proof of Theorem 4.5]{Alexeev-1994ge} and so 
$K_{\tilde X}$ is Cartier. Let $\tilde f: \tilde X\to \bar X\subset \PP^n$
be the morphism defined by $|-K_{\tilde X}|$ and let $\tilde X \to Y \to \bar X$ be the Stein factorization.
We have a $G$-equivariant diagram
\[
\xymatrix{
&\tilde X\ar[dl]_{f}\ar[dr]^{h}\ar@/^0.9pc/[drrr]^{\tilde f}&
\\
X\ar@{-->}[rr]&&Y\ar[rr]^{} && \bar X\subset \PP^n
} 
\]
Let $\tilde H\in |-K_{\tilde X}|$ be a general member.
Then $\tilde H$ is a smooth K3 surface.
Since the restriction map 
$H^0(\tilde X,-K_{\tilde X})\to H^0(\tilde H,-K_{\tilde X})$ is surjective,
by \cite{Saint-Donat-1974} 
the image $\tilde f(\tilde H)\subset \PP^n$ is either 
\begin{enumerate}
\renewcommand\labelenumi{(\arabic{enumi})}
\renewcommand\theenumi{(\arabic{enumi})}
\item 
a rational normal curve,
\item 
a rational normal surface of degree $n-1$ (see \ref{Theorem-Del-Pezzo-Bertini}), or 
\item 
a K3 surface with Du Val singularities. 
\end{enumerate}

Consider the case $\dim \bar X=1$. 
Then $\bar X=\tilde f(\tilde H)$ is a rational normal curve and 
$-K_{\tilde X}=\tilde f^*\OOO_{\bar X}(1)=n\tilde f^*\tilde P$,
where $\tilde P\in \bar X$ is a point. By Corollary \ref{Corollary-M1}
$n=1$ and we get (i).

Consider the case $\dim \bar X=2$. 
Then we are in situations (1) or (2) above.
Both possibilities give us that 
$\bar X\subset \PP^n$ is a rational normal surface of degree $n-1$.
Since
$-K_{\tilde X}=\tilde f^*\OOO_{\bar X}(1)$, again by Corollary \ref{Corollary-M1} the linear system
$|\OOO_{\bar X}(1)|$ cannot be decomposed in a sum of movable linear systems.
The only possibility is (ii) (see Theorem \ref{Theorem-Del-Pezzo-Bertini}). 

Finally, in the case $\dim \bar X=3$ we get (iii) as in \cite[Corollary 4.6]{Alexeev-1994ge}. 
\end{proof}

\begin{scase}{\bf Example.}
Consider famous Fletcher's list of 95 terminal Fano weighted hypersurfaces \cite{Iano-Fletcher2000}.
92 of them satisfy the condition $\dim |-K_X|\le 2$.
Furthermore, many of these varieties admit a biregular involution $\tau$
with $\kappa(\tau,X)\ge 0$.

For example a hypersurface 
$X=X_{66}\subset\PP(1,5,6,22,33)$
in a suitable coordinate system $(x_1,x_5,x_6,x_{22},x_{33})$ 
can be given by the equation
$x_{33}^2=\phi(x_1,x_5,x_6,x_{22})$, where $\phi$ is a quasihomogeneous 
polynomial of weighted degree $66$.
The projection $X\to \PP(1,5,6,22)$ is a double cover 
and for the corresponding Galois involution $\tau$ we have 
$\Fix(\tau)\simeq \{\phi=0\}\subset \PP(1,5,6,22)$.
It is easy to check that for a general choice of 
$\phi$ the subvariety $S:=\Fix(\tau)$ is a normal surface with 
only cyclic quotient singularities of types $\frac 12(1,1)$ and $\frac15(1,2)$.
By the adjunction formula, $K_S=\OOO_{S}(32)$ and so $K_S^2=512/5$.
Let $\mu: \tilde S\to S$ be the minimal resolution.
Then $K_{\tilde S}=\mu^*K_S-\Delta$, where $\Delta$ is the codiscrepancy divisor.
In our case $\Delta$ is supported over the point of type $\frac15(1,2)$,
so $\Delta ^2=-2/5$ and $K_{\tilde S}^2=512/5-2/5=102$.
Therefore, $S$ is a surface of general type and $\kappa(\tau,X)=2$. 

Note however that a general Fano hypersurface $X$ in the Fletcher's list
is not rational \cite{Corti2000}. So, for general $X$ our construction does not
produce any involution in the Cremona group.
It is interesting to investigate \textit{special} members of this list.
\end{scase}

\section{Gorenstein Fano threefolds: examples}
\label{section-Gorenstein-examples}
In this section we collect several examples of
involutions acting on Gorenstein canonical Fano threefolds.
We are interested only in those involutions 
that are not conjugate to actions on conic bundles nor del Pezzo fibrations.
First of all, a hyperelliptic Fano threefold has the Galois 
involution:

\begin{case}{ \bf Example.} \label{example-hyperelliptic-double-solid}
Let $\pi: X\to \PP^3$ be 
a double cover branched over a divisor $B\subset \PP^3$
of degree $6$ such that the pair $(\PP^3,\frac12B)$ is klt.
Then $X$ is a hyperelliptic Fano threefold with $\g(X)=2$.
Let $\tau$ be the Galois involution.
Then $\Sing(X)=\pi^{-1}(\Sing(B))$ and
$\Fix(\tau)=\pi^{-1}(B)$. If $B$ is smooth or has only Du Val singularities,
then $\kappa(\tau,X)=2$.
Note however that if $X$ has only isolated cDV singularities 
and is $\QQ$-factorial, it cannot be rational \cite{Cheltsov2010a}.
In particular, if we are interested in elements of order $2$ in $\Cr_3(\Bbbk)$,
we should consider singular sextics $B\subset \PP^3$.

For example, let $B\subset \PP^3$ be the \textit{Barth sextic}.
Recall that it is a surface given by the equation
\[
4(\epsilon^2 x_1^2-x_2^2)(\epsilon^2 x_2^2-x_3^2)(\epsilon^2 x_3^2-x_1^2)-x_4^2 (1+2\epsilon)(x_1^2+x_2^2+x_3^2-x_4^2)^2=0,
\] 
where $\epsilon=(1+2 \sqrt 5)/2$. Then $X$ has 
exactly $65$ nodes and has no other singularities. Moreover, 
$X$ is rational (see \cite[Example 3.7]{Endrass1999}). This gives us an example of involution $\tau\in\Cr_3(\Bbbk)$
with $\kappa(\tau,X)=2$.
\end{case}

\begin{case}{\bf Example.} \label{example-hyperelliptic-double-quadric}
Let $W=W_2\subset \PP^4$ be a smooth quadric and let $\pi: X\to W$ be 
a double cover branched over a divisor $B\in |\OOO_W(4)|$
such that the pair $(W,\frac12B)$ is klt.
Then $X$ is a hyperelliptic Fano threefold with $\g(X)=3$.
Let $\tau \in \Aut(X)$ be the Galois involution.
As above, $\kappa(\tau,X)=2$ if the branch divisor has mild singularities, however, 
$X$ is not rational if it is smooth \cite{Iskovskikh1980}. 
\end{case}

A series of examples comes from double covers of del Pezzo threefolds:

\begin{case}{ \bf Construction.} \label{Construction-del-Pezzo-threefolds}
Let $Y$ be a del Pezzo threefold with only canonical
Gorenstein singularities. By definition, 
$-K_Y\sim 2A$ for some $A\in \Pic(Y)$.
Let $B\in |-K_Y|$ be an element such that the pair 
$(Y,\frac 12 B)$ is klt. 
Consider the double cover $\pi: X\to Y$ branched over $B$.
By the Hurwitz formula 
\[
-K_X=-\pi^*\left(K_Y+\frac 12 B\right)=\pi^*A.
\]
Then $X$ is a Fano threefold with 
canonical Gorenstein singularities.
Denote by $\tau$ the Galois involution and $S:=\pi^{-1}(B)$.
Then $S\simeq B$, $\Sing(X)=\Sing(S)\cup \pi^{-1}\Sing(Y)$, and
$\Fix(\tau)=S$. Since $K_B=0$, we have $\kappa(\tau,X)\le 0$ and 
$\kappa(\tau,X)= 0$ if and only if $B$ is irreducible, normal 
and has only Du Val singularities. 
Let $d=A^3$ be the degree of $Y$. 
We have $-K_X^3=2A^3$, so $\g(X)=d+1$.
\end{case}

Consider particular cases of our construction.

\begin{scase}{\bf Example ($g=2$).}\label{example-hyperelliptic-double-space-linear}
Let $Y$ be a del Pezzo threefold of degree $1$.
Then $Y$ is a hypersurface of degree $6$ in $\PP(1^3,2,3)$.
The projection $\delta: Y\to \PP(1^3,2)=:\PP$ is a double cover.
Consider another double cover $\gamma: \PP^3\to \PP$, the quotient by
a reflection. Let $X:=Y\times _{\PP} \PP^3$. Our involution $\tau$ is
the Galois involution of the projection $X\to Y$ which is a double cover.
In other words, $X$ is a hypersurface of degree $6$ in $\PP(1^4,3)$
given by the (quasihomogeneous) equation $z^2=\phi(x_0^2,x_1,x_2,x_3)$ and $\tau$ is the reflection 
$x_0\mapsto -x_0$. 
\end{scase}

\begin{scase}{\bf Example ($g=3$).}\label{example-hyperelliptic-double-quadric-linear}
Let $Y$ be a del Pezzo threefold of degree $2$.
Then $Y$ is a hypersurface of degree $4$ in $\PP(1^4,2)$.
The projection $\delta: Y\to \PP(1^4)=\PP^3$ is a double cover.
As above, consider another double cover $\gamma: W_2\to \PP$, the projection of a quadric $W_2\subset \PP^4$ 
away from a point $P\notin W_2$. Let $X:=Y\times _{\PP^3} W_2$. Our involution $\tau$ is
the Galois involution of the projection $X\to Y$ which is a double cover.
\end{scase}

\begin{scase}{\bf Example ($g=3$).}\label{example-hyperelliptic-quartic}
Let $X\subset \PP^4$ be a quartic given by the (homogeneous) equation 
$ax_0^4+x_0^2\phi_2(x_1,x_2,x_3,x_4)+\phi_4(x_1,x_2,x_3,x_4)=0$.
Our involution $\tau$ is the reflection 
$x_0\mapsto -x_0$. 
The quotient is given by 
$ay^2+y\phi_2(x_1,x_2,x_3,x_4)+\phi_4(x_1,x_2,x_3,x_4)=0$ in $\PP(1^4,2)$.
If $a\neq 0$, this is a del Pezzo threefold of degree $2$.
If $a=0$, then $\tau$ is linearizable.
\end{scase}

\begin{scase}{\bf Example ($g=4$).}\label{example-hyperelliptic-V6}
Let $X=X_{2\cdot 3}\subset \PP^5$
be an intersection of a quadric $Q$ and a cubic cone $V$
with vertex $P\in V$ so that $P\notin Q$.
The projection $\pi: X\to Y\subset \PP^4$ away from $P$ is a double cover
of a cubic $Y\subset \PP^4$. 
Then $\tau$ is
the Galois involution of this projection.
\end{scase}

\begin{scase}{\bf Example ($g=5$).}\label{example-hyperelliptic-V8}
Let $X=X_{2\cdot 2\cdot 2}\subset \PP^6$
be an intersection of three quadrics $Q_1$, $Q_2$, $Q_3$ 
so that $Q_1$ and $Q_2$ are cones with vertices at the same point 
$P\in Q_1\cap Q_2$, where $P\notin Q_3$.
The projection $\pi: X\to Y\subset \PP^5$ away from $P$ is a double cover
of an intersection of two quadrics $Y\subset \PP^5$. 
As above $\tau$ is
the Galois involution of this projection.
\end{scase}

\section{Gorenstein Fano threefolds: the proof of the main theorem}
\label{section-Gorenstein} 
\begin{case}{\bf Assumptions.}\label{Notation-Gorenstein}
Throughout this section $X$ is a Fano threefold with 
only canonical Gorenstein singularities. 
Let $g:=\g(X)\ge 2$ be the genus of $X$. So $-K_X^3=2g-2$ and $\dim |-K_X| = g+1$.
Let $\Phi=\Phi_{|-K_X|}: X \dashrightarrow \PP^{g+1}$ be the anti-canonical map.
Assume that $(X,G)$ is not of type \type{C} nor \type{D}
(otherwise we are in the situation of \S \ref{section-Conic-bundles} or
\S \ref{section-Del_Pezzo}).
Our aim is to replace $(X,G)$ with another birational model satisfying 
the above assumptions and \textit{minimal} possible genus. 
\end{case}

\begin{mtheorem} {\bf Lemma.}\label{Lemma-cDV}
The singularities of $X$ are of type cDV.
\end{mtheorem}
\begin{proof}
Recall that a three-dimensional canonical Gorenstein singularity 
is of type cDV if and only if the center of every crepant exceptional divisor 
is one-dimensional.
Let $P\in X$ be a non-cDV singularity. For its orbit $\Lambda$
we have two possibilities: $\{P\}$ and $\{P,\, P'\}$, where $P'\neq P$.
Take the subsystem $\LLL \subset |-K_X|$ consisting of all the members passing through $\Lambda$.
Since $\dim |-K_X| = g+1\ge 3$, $\dim \LLL\ge 1$. 
Write $\LLL=D+\MMM$, where 
$D$  (resp. $\MMM$) is the fixed (resp.  movable) part of $\LLL$.
By our assumption $(X,G)$ is not of type \type{C} nor \type{D}
and by Lemma \ref{Lemma-index>1} we have $D=0$ and $\MMM=\LLL$.
Let $E$ be a crepant exceptional divisor with center at $P$.
Since $P\in \Bs \LLL$, for the discrepancy of $E$ we have 
$a(E,X,\LLL)< a(E,X,0)=0$. Hence, the pair $(X,\LLL)$ is not canonical.
This contradicts  Lemma \ref{Lemma-index>1}.
\end{proof}

\begin{mtheorem} {\bf Proposition.}\label{Proposition-hyperelliptic}
Assume that $|-K_X|$ is not very ample. 
Then $X$ is one of the following:
\begin{enumerate}
\item
$X$ is a double cover of $\PP^3$ branched over a sextic surface $B\subset \PP^3$;
\item
$X$ is a double cover of a \emph{smooth} quadric $W_2\subset \PP^4$
branched over a surface cut out on $X$ by a quartic hypersurface.
\end{enumerate}
\end{mtheorem}

\begin{proof}
First we assume that $\Bs |-K_X|\neq \emptyset$. Then the 
image $W:=\Phi(X)\subset \PP^{g+1}$ is a surface of minimal degree (see \ref{Fano-threefolds}).
Apply Theorem \ref{Theorem-Del-Pezzo-Bertini}.
Since $\dim |-K_X|=g+1\ge 3$, the surface $W$ is either a Hirzebruch surface $\FF_n$ 
or a cone $\PP(1,1,n)$ over a rational normal curve.
In both cases we get a contradiction by Corollary \ref{Corollary-M1}
because $|-K_X|=\Phi^* |\OOO_W(1)|$ is decomposable in a sum of two movable linear systems.

Thus $X$ is hyperelliptic and the linear system $|-K_X|$ defines 
a double cover $\Phi :X\to W$, where $W\subset \PP^{g+1}$ is a (normal) variety of degree $g-1$ 
(see \ref{Fano-threefolds}).
Apply Theorem \ref{Theorem-Del-Pezzo-Bertini}.
In the case (3) $W$ is a Veronese cone, $W\simeq \PP(1^3,2)$.
Thus $-K_{X}= 2\Phi^* \OOO_{\PP(1^3,2)} (1)$, where $|\OOO_{\PP(1^3,2)} (1)|$
is a movable linear system. This contradicts Corollary \ref{Corollary-M1}.
In the case (4), 
we consider the following diagram:
\[
\xymatrix {
X \ar[d] & \tilde X\ar[l]\ar[d]\ar[dr]
\\
W & \PP(\EEE)\ar[l]\ar[r]& \PP^1
} 
\]
where $\tilde X$ is the normalization of the product $X\times_W \PP(\EEE)$.
Then $\tilde X$ is a weak Fano threefold with canonical Gorenstein singularities 
\cite[Lemma 3.6]{Przhiyalkovskij-Chel'tsov-Shramov-2005en}
and the morphism $\tilde X\to X$ is crepant and $G$-equivariant.
Hence the composition $\tilde X\to \PP^1$ is a $G$-equivariant del Pezzo
fibration (of degree $2$). This contradicts our assumptions. 
In the case (1) we get (i)  (see \cite{Przhiyalkovskij-Chel'tsov-Shramov-2005en}).
Finally, in the case (2) we have a double cover as in (ii).
We have to show only that the quadric $W_2\subset \PP^4$ is smooth.
Indeed, otherwise 
the linear system of its hyperplane sections $|\OOO_{W_2}(1)|$ admits a decomposition in 
a sum of two $G$-invariant movable linear systems of Weil divisors 
and so $|-K_X|=\Phi^*|\OOO_{W_2}(1)|$ is.
This again contradicts Corollary \ref{Corollary-M1}.
\end{proof}

\begin{mtheorem} {\bf Proposition.}\label{Proposition-quadric}
Assume that $X$ is trigonal. 
Then $X$ is one of the following:
\begin{enumerate}
\item
$Y\simeq Y_{4} \subset \PP^4$;
\item
$Y\simeq Y_{2\cdot 3} \subset \PP^5$.
\end{enumerate}
\end{mtheorem}

\begin{proof}
We have a $G$-equivariant
anti-canonical embedding $X=X_{2g-2}\subset \PP^{g+1}$. Assume that 
$X$ is not an intersection of quadrics. Then the quadrics through $X$ in $\PP^{g+1}$ 
cut out a fourfold $W\subset\PP^{g+1}$ of degree $g- 2$ (see \ref{Fano-threefolds}). As above,
apply Theorem \ref{Theorem-Del-Pezzo-Bertini}.
If $W=\PP^4$ or $W=W_2\subset \PP^5$, then we get (iii) and (iv) respectively.
In the remaining cases, let $\sigma: \tilde W\to W$ be the blow up 
of the singular locus (we put $\sigma=\operatorname{id}$ 
if $W$ is smooth) and let $\tilde X\subset \tilde W$ be the proper transform of $X$.
Then $\tilde W$ is smooth, 
$\tilde X$ is a weak Fano threefold with canonical Gorenstein singularities 
\cite[Lemmas 4.4, 4.7]{Przhiyalkovskij-Chel'tsov-Shramov-2005en},
and the morphism $\tilde X\to X$ is crepant (and $G$-equivariant).
The variety $\tilde W$ is either $\PP^2$-bundle (resp. $\PP^3$-bundle)
over $\PP^2$ (resp $\PP^1$). By the adjunction formula, the projection to the base induces 
on $\tilde X$ a $G$-equivariant structure of a
conic bundle (resp. del Pezzo fibration). This contradicts our assumptions. 
\end{proof}

\begin{mtheorem} {\bf Lemma.}\label{Lemma-lines}
Assume that $|-K_X|$ is very ample. Then 
there are no lines passing through a general point $P\in X$.
\end{mtheorem}
\begin{proof}
Indeed, otherwise the family of lines on
$X$ has dimension at least $2$. Let $L\subset X$ be a line passing through $P$.
Let $\delta: X'\to X$ be a terminal crepant model and let $L'\subset X'$ be the proper transform of $L$.
If $L'\cap \Sing(X')\neq \emptyset$, then 
since $\dim \Sing(X')\le 0$,
we see that $X$ is a cone over $S$. 
The base of this cone is a hyperplane section, a K3 surface
with only Du Val singularities. 
But then the vertex is not a canonical singularity, 
a contradiction.
Hence $L'$ is contained into the smooth locus of $X'$.
Since $L'$ belongs to a covering family of rational curves,
we have $N_{L'/X'}\simeq \OOO_{\PP^1}(a) \oplus \OOO_{\PP^1}(b)$, where $a,b \ge 0$. 
On the other hand, $\deg N_{L'/X'}= -K_{X'}\cdot L'-2=-1$, a contradiction.
\end{proof}

\begin{mtheorem} {\bf Lemma.}\label{Lemma-projection-1}
Assume that $\g(X)\ge 3$ and there exists a $G$-fixed singular point $P\in X$.
Then there exists a $G$-equivariant birational map
$X \dashrightarrow Y$, where $Y$ is a Fano threefold with canonical Gorenstein singularities
and $\g(Y)< \g(X)$. 
\end{mtheorem}
\begin{proof}
By Proposition \ref{Proposition-hyperelliptic} either 
$-K_X$ is very ample or $X$ is of type \ref{Proposition-hyperelliptic}(ii),
i.e. $X$ is hyperelliptic and $\Phi(X)$ is a smooth quadric
$W_2\subset \PP^4$. In the second case $P$ lies on the ramification divisor
because $\Phi$ cannot be \'etale over $\Phi(P)$.

Let $\MMM\subset |-K_X|$ be the subsystem consisting of all the members passing through $P$.
Thus $\Bs |\MMM|=\{P\}$. 
We claim that the map $h: X \dashrightarrow \bar X\subset \PP^g$ 
given by $\MMM$ is generically finite.
If $-K_X$ is very ample, then $h$ is nothing but the projection away from $P$.
Since $X$ is not a cone, $h$ must be generically finite.
So we may assume that 
$\MMM=\Phi^*\LLL$, where $\LLL$ is the linear system of
hyperplane sections of $W_2$ passing through $P':=\Phi(P)$. We have the following commutative diagram:
\[
\xymatrix{
X\ar[d]^{\Phi}\ar@{-->}[rr]^{h}&& \bar X\ar[d]
\\
W_2\ar@{-->}[rr]^{}&& \PP^3
} 
\]
where $W_2\dashrightarrow \PP^3$ is the projection away from $P'$.
Hence $h$ is generically finite.

If the pair $(X,\MMM)$ is not canonical, then by Lemma 
\ref{Lemma-index>1} $(X,G)$ is of type \type{C} or \type{D}.
Thus $(X,\MMM)$ is canonical (but not terminal at $P$).
Consider its log crepant terminal model 
$f:(\tilde X,\tilde \MMM)\to (X,\MMM)$ (see \ref{pairs}).
Thus $(\tilde X,\tilde \MMM)$ is terminal, 
$$
K_{\tilde X}+\tilde \MMM=f^*(K_X+\MMM)\sim 0,
$$
and so $\tilde \MMM\subset |-K_{\tilde X}|$.
Write
$$
K_{\tilde X}=f^*K_X+\Theta,\qquad \tilde \MMM=f^*\MMM-\Theta,
$$
where $\Theta$ is an effective exceptional divisor
such that $\Supp(\Theta)=f^{-1}(P)$.
For a general member $D\in |-K_{\tilde X}|$ we have 
$f_*D\in |-K_X|= |-f^*K_X-\Theta|$. Since $\Theta$ cannot be numerically trivial over $X$,
we have $\Supp(\Theta)\cap \Supp(D)\neq \emptyset$,
$P\in f(D)$, and $f_*D\in \MMM$.
Therefore, $\tilde \MMM= |-K_{\tilde X}|$.

Since $(\tilde X,\tilde \MMM)$ is terminal, 
the base locus of $\tilde \MMM$ consists of at most
a finite number of points   \cite[Lemma 1.22]{Alexeev-1994ge}.
In particular, $-K_{\tilde X}$ is a nef Cartier divisor. By the above $-K_{\tilde X}$ is
also big, i.e. $\tilde X$ is a weak Fano threefold.
Consider the anti-canonical model $Y:=\Proj \oplus_{n\ge 0} H^0(\tilde X, -nK_{\tilde X})$.
Then $Y$ is a Fano threefold with only canonical Gorenstein singularities.
Let $\tilde X\to Y \dashrightarrow \bar X$ be the Stein factorization
of  $h\circ f$.
Thus we have a diagram
\begin{equation}\label{equation-Lemma-projection-1}
\xymatrix{
&\tilde X\ar[dl]_{f}\ar[dr]^{}\ar@{-->}@/^0.9pc/[drrr]^{}&
\\
h:X\phantom{h}\ar@{-->}[rr]^{}&&Y\ar@{-->}[rr]^{} && \bar X
} 
\end{equation}
By Proposition \ref{Proposition-hyperelliptic} the linear system 
$|-K_Y|$ is base point free and so $|-K_{\tilde X}|$ is.
Our construction is $G$-equivariant. 
Finally,
\[
\dim |-K_Y|=\dim |-K_{\tilde X}|=\dim \tilde \MMM=\dim \MMM=
\dim |-K_{X}|-1. 
\]
Therefore, $\g(Y)<\g(X)$. 
\end{proof}

\begin{mtheorem} {\bf Lemma.}\label{Lemma-projection-2}
Assume that $\g(X)\ge 4$ and $X$ is singular.
Then there exists a $\tau$-equivariant birational map
$X \dashrightarrow Y$, where $Y$ is a Fano threefold with only canonical Gorenstein singularities
and $\g(Y)< \g(X)$. 
\end{mtheorem}
\begin{proof}
By Proposition \ref{Proposition-hyperelliptic} the linear system $|-K_X|$ 
defines an embedding $X=X_{2g-2}\subset \PP^{g+1}$. 
By Lemma \ref{Lemma-projection-1} we may assume that $\tau$ has no fixed singular points.
So, there are two points $P_1$, $P_2$ that are switched by $\tau$.
Then, similar to the proof of Lemma \ref{Lemma-projection-1},
we consider the subsystem $\MMM\subset |-K_X|$ consisting of all the members passing through $P_1$ and $P_2$.
It gives a map $h: X\dashrightarrow \bar X\subset \PP^{g-1}$ which  is 
the projection away from the line $L\subset \PP^{g+1}$ passing through $P_1$ and $P_2$.
We claim that $h$ is generically finite.
Assume that its general fiber contains a curve $C$.
Then $C\subset \Pi \cap X$, where $\Pi$ is a plane passing through $L$.
Since $(X,G)$ is not of type \type{C} the curve $C$ is not rational.
In this case, $X$ cannot be an intersection of quadrics. By Proposition \ref{Proposition-quadric}
the only possibility is that $X$ is an intersection of a 
quadric $Q$ and a cubic $Y$ in $\PP^5$. Again since $C\subset \Pi \cap X$ is a non-rational curve,
we have $\Pi\subset Q$. If $L\not \subset \Sing(Q)$, there exists only a finite number 
of planes in $Q$ passing through $L$. Thus we may assume that $L \subset\Sing(Q)$.
But then $L\cap Y$ either coincides with $L$ or 
is an invariant $0$-cycle of degree $3$.
In both cases $\tau$ has a fixed point on $L$ which must be a singular 
point of $X$. This contradicts out assumptions.

Therefore, $h$ is generically finite.
Then, similar to the proof of Lemma \ref{Lemma-projection-1},
we get the diagram \eqref{equation-Lemma-projection-1}
and we can show that $\dim |-K_Y|<\dim |-K_X|$.
\end{proof}

\begin{mtheorem} {\bf Lemma.}\label{Lemma-hyperplane}
Assume that $S:=\F(\tau)\neq \emptyset$.
Furthermore, assume that $\g(X)$ is minimal among all birational models of
$(X,G)$ satisfying \xref{Notation-Gorenstein}.
Then either 
\begin{enumerate}
 \item 
$X$ is hyperelliptic and
$\tau$ is the hyperelliptic involution or 
\item
$S\in |-K_X|$ and $S$ has Du Val singularities only.
\end{enumerate}
\end{mtheorem}

\begin{proof}
By Proposition \ref{Proposition-hyperelliptic} we may assume that 
$|-K_X|$ defines a morphism $\Phi: X\to \PP^{g+1}$.
If the natural action of $\tau$ on
$H^0(X,-K_X)$ is trivial, 
then $X$ must be hyperelliptic and we get the case (i).

Thus we may assume that $\tau$ faithfully acts on
$H^0(X,-K_X)$. Then $\Phi(S)$ is contained in a hyperplane
$\PP^g\subset \PP^{g+1}$.
So, we can write 
$\Phi^{-1}(\PP^g)=S+R\sim -K_X$, where $R$ is an effective Weil divisor.
Let $\nu: S'\to S$ be the normalization.
By the subadjunction formula (see e.g. \cite[Ch. 16]{Utah})
\[
K_{S'}+\Diff_{S'} (R) =\nu^* (K_X+S+R)=0, 
\]
where $\Diff$ is the so-called \emph{different}, an effective $\QQ$-divisor.
By \cite[Cor. 3]{Miyaoka1986} we have $\Diff_{S'} (R)=0$.
So $\dim S\cap R\le 0$ and $S$ is smooth in codimension one.
Since $S+R$ is Cohen-Macaulay, we get that $R=0$ and $S$ is normal.
Similarly, considering the minimal resolution of $S$, 
we obtain that the singularities of  $S$ are
Du Val. 
\end{proof}

\begin{mtheorem} {\bf Lemma.}\label{Lemma-projection-3}
Assume that $\g(X)\ge 6$ and $\F(\tau)\neq\emptyset$.
Then
there exists a $G$-equivariant birational map
$X \dashrightarrow Y$, where $Y$ is a Fano threefold with canonical Gorenstein singularities
and $\g(Y)< \g(X)$. 
\end{mtheorem}

\begin{proof}
By Lemma \ref{Lemma-lines} the union of all the lines on $X$
is a closed proper subset $V$.
Since the surface $\F(\tau)$ is not uniruled, 
it is not a component of $V$.
Take a point $P\in \F(\tau)\setminus V$.
Let $f: \tilde X\to X$ be the blowup of $P$ and let 
$E:=f^{-1}(P)$ be the exceptional divisor.
We claim that the divisor $-K_{\tilde X}=-f^*K_X-2E$
is nef. Indeed, if $-K_{\tilde X}\cdot \tilde C<0$ for some curve $\tilde C$,
then $\tilde C$ is contained in $\Bs |-K_{\tilde X}|$. 
Note that $\MMM:=f_*|-K_{\tilde X}| \subset |-K_X|$ 
is the subsystem consisting of
all the members singular at $P$. 
Hence  for 
$C:=f(\tilde C)$ we have $C\subset X\cap T_{P,X}$. By 
Proposition \ref{Proposition-quadric}
the variety $X$ is an intersection of quadrics.
Hence $C$ is a line passing through $P$.
This contradicts our choice of $P$.
Thus $-K_{\tilde X}$ is nef.
Since $-K_{\tilde X}^3= -K_X^3-8>0$, it is big.
By the base point free theorem (see e.g. \cite[Th. 3.3]{Kollar-Mori-19988}) the linear system 
$|-nK_{\tilde X}|$ defines a morphism $h:\tilde X\to Y\subset \PP^n$
and we may assume that $Y$ is normal, has only canonical Gorenstein singularities,
$K_{\tilde X}=h^*K_Y$,
and $-K_Y$ is ample.
Thus we get  
a new Fano threefold $Y$ with $-K_{Y}^3= -K_X^3-8>0$.
\end{proof}

From the above assertions we get a description of our Fano 
threefold $X$.

\begin{mtheorem} {\bf Corollary.}\label{Corollary-Gorenstein-classification}
Replacing $(X,\tau)$ with a birational model we may assume that $\g(X)\le 5$ and one of the following holds:
\begin{enumerate}
\item
$X$ is a double cover of $\PP^3$ 
branched over a sextic $B\subset \PP^3$, the singularities of $X$ are cDV;

\item
$X$ is a double cover $X$ of a smooth quadric $W_2\subset \PP^4$ 
branched over a surface $B\subset W_2$ of degree $8$, the singularities of $X$
are terminal;

\item
$X\subset \PP^4$ is a quartic with terminal singularities;
\item
$X\subset \PP^5$ is a smooth intersection of a quadric and a cubic cone;
\item
$X\subset \PP^6$ is a smooth intersection of three quadrics.
\end{enumerate}
\end{mtheorem}
\begin{proof}
By Lemma \ref{Lemma-cDV} the singularities of $X$ are cDV.
Recall that $\g(X)\ge 2$ (because $-K_X^3= 2\g(X)-2>0$)
and if $\g(X)= 2$, then $\dim |-K_X|=3$ (see \ref{Notation-Gorenstein})
and we get case (i) by Proposition \ref{Proposition-hyperelliptic}.
Thus we assume that $\g(X)\ge 3$.
By the classification of 
Applying Lemma \ref{Lemma-projection-3} we get a model $(X,G)$ with $\g(X)\le 5$.
In the cases $\g(X)=5$ and $4$ by Lemma \ref{Lemma-projection-2}
we may assume that $X$ is smooth.
In the case $\g(X)=3$ by Lemma \ref{Lemma-projection-1}
we may assume that $\Fix(\tau,X)\cap \Sing(X)=\emptyset$,
in particular, $X$ has only isolated cDV singularities.
\end{proof}

In the case \type{F^c} (a) we have an embedding 
$X \hookrightarrow \PP(1,1,1,1,3)$ 
which is equivariant (because this embedding is given by
the sections of $H^0(X,-nK_X)$). In suitable coordinates 
we get the desired subcases (i) and (ii).
It remains to show that in the cases \type{F^c} (b)(ii), (c), and (d)
the action of $\tau$ is given by Construction \ref{Construction-del-Pezzo-threefolds}.
This is a consequence of the following.

\begin{mtheorem} {\bf Lemma.}
Let $X$ be as in 
\textup{(ii)-(v)} of \xref{Corollary-Gorenstein-classification}.
Assume that $\F(\tau,X)\neq \emptyset$ and the action of $G$ on
$H^0(X,-K_X)$ is not trivial. 
Furthermore, assume that $\g(X)$ is minimal among all birational models of
$(X,G)$ satisfying \xref{Notation-Gorenstein}.
Then the quotient $X/G$ is a \textup(Gorenstein\textup) del Pezzo threefold
of degree $\g(X)-1$.
\end{mtheorem}

\begin{proof} 
Let $\pi: X\to Y=X/G$ be the quotient map,
let $S:=\F(\tau, X)$,
and let $R:=\pi(S)$.
By Lemma \ref{Lemma-hyperplane} we have $S\sim -K_X$. 
 The divisor $-K_X+S=-\pi^*K_Y$ is ample.
Hence, $Y$ is a (log terminal) Fano threefold.
By the adjunction formula, $0=K_R=(K_Y+R)|_R$. Hence, $K_Y+R\qq 0$
and $-K_Y^3=4 (-K_X)^3=8(g-1)$. 
Moreover, the branch divisor $R\subset Y$ is divisible by $2$ in $\Cl(Y)$. 
Thus it is sufficient to show that $Y$ has only 
terminal Gorenstein singularities.

By Lemma \ref{Lemma-projection-1} $\Fix(\tau,X)\cap \Sing(X)=\emptyset$.
Therefore, $S$ is smooth and so
$Y$ is smooth along $R$.
Assume that $\Fix(\tau,X)\neq S$. Then $\Fix(\tau,X)\setminus S$ consists of a finite number of smooth points
$P_1,\dots,P_k\in X$ and by the above each $\pi(P_i)\in Y$ is a point of type $\frac12(1,1,1)$. 
Then the linear span of $\Phi(S)$ is a hyperplane $\PP^g$
and $\Phi(P_1)=\cdots=\Phi(P_k):=Q\notin \Phi(S)$.
Here $\Fix(\tau,\PP^{g+1})= \PP^g\cup \{Q\}$.
In particular, $k\le 2$.
Easy computations with the orbifold Riemann-Roch formula
show that the group $\Cl(Y)$ is torsion free (see e.g. \cite[Prop. 2.9]{Prokhorov2008a}). 
Since $K_Y+R\qq 0$, the divisor $K_Y$ is Cartier, a contradiction.
Thus $\Fix(\tau,X)= S$ and the singularities of $Y$ are terminal Gorenstein.
\end{proof}

\par \smallskip\noindent
{\bf Acknowledgments.}
The author would like to thank the referee for numerous 
helpful comments.

 \end{document}